# The Computation of Key Properties of Markov Chains via Perturbations

Jeffrey J. Hunter

*Department Mathematical Sciences*
*School of Engineering, Computer and Mathematical Sciences, Auckland University of Technology,*
*Private Bag 92006, Auckland 1142, New Zealand*

## Abstract

Computational procedures for the stationary probability distribution, the group inverse of the Markovian kernel and the mean first passage times of a finite irreducible Markov chain, are developed using perturbations. The derivation of these expressions involves the solution of systems of linear equations and, structurally, inevitably the inverses of matrices. By using a perturbation technique, starting from a simple base where no such derivations are formally required, we update a sequence of matrices, formed by linking the solution procedures via generalized matrix inverses and utilising matrix and vector multiplications. Four different algorithms are given, some modifications are discussed, and numerical comparisons made using a test example. The derivations are based upon the ideas outlined in Hunter, J.J., "The computation of stationary distributions of Markov chains through perturbations", Journal of Applied Mathematics and Stochastic Analysis, 4, 29-46, (1991).

## 1. Introduction

In Markov chain theory stationary distributions, mean first passage times and the group inverse provide significant information regarding the behaviour of the chain.

Let $\{X_n, n \geq 0\}$ be a finite Markov chain (M. C.) with state space $S = \{1, 2, \ldots, m\}$ and transition matrix $P = [p_{ij}]$, where $p_{ij} = P\{X_n = j \mid X_{n-1} = i\}$ for all $i, j \in S$.

It is well known (Feller (1950), Kemeny and Snell (1960), that if the M. C. is *regular* (irreducible and aperiodic) then for all $i$, $j$, $\lim_{n\to\infty} p_{ij}^{(n)} = \lim_{n\to\infty} p_j^{(n)} = \pi_j$ where $p_{ij}^{(n)} = P\{X_n = j \mid X_0 = i\}$, $p_j^{(n)} = P\{X_n = j\}$. The limiting probability of being in state $j, \pi_j$, is in fact the "*stationary probability*" of being in state $j$, in that if $P\{X_0 = j\} = \pi_j$ for all $j$, then $P\{X_n = j\} = \pi_j$, for all $j$ and $n \geq 0$. An important result is that the stationary distribution $\{\pi_j\}$, $(1 \leq j \leq m)$, exists and is unique for all irreducible M. C.'s, that $\pi_j > 0$ for all $j$, and satisfies the equations (the *stationary equations*)

$$\pi_j = \sum_{i=1}^{m} \pi_i p_{ij} \text{ with } \sum_{i=1}^{m} \pi_j = 1. \tag{1.1}$$

If $\boldsymbol{\pi}^T = (\pi_1, \pi_2, \ldots, \pi_m)$, the stationary probability vector, and $\boldsymbol{e}$ is a column vector of 1's, the stationary equations (1.1) can be expressed as

$$\boldsymbol{\pi}^T (I - P) = \boldsymbol{0}^T, \text{ with } \boldsymbol{\pi}^T \boldsymbol{e} = 1. \tag{1.2}$$

Thus $\boldsymbol{\pi}^T$ can be determined by solving a constrained system of linear equations involving the singular matrix $I - P$ (since each row of $P$ is a discrete distribution, and $P$ is a stochastic matrix with each row sum 1, i.e. $P\boldsymbol{e} = \boldsymbol{e}$).

Let $\Pi = \boldsymbol{e}\boldsymbol{\pi}^T$. In the case of a regular M. C. (finite, irreducible and aperiodic), $\lim_{n\to\infty} P^n = \Pi$, and, in the case of a finite irreducible M. C., $\lim_{n\to\infty} \dfrac{I + P + P^2 + \ldots + P^n}{n} = \Pi$.

Let $T_{ij} = \min\left[n \geq 1, X_n = j \mid X_0 = i\right]$ be the first passage time from state $i$ to state $j$ (first return when $i = j$) and define $m_{ij} = E\left[T_{ij} \mid X_0 = i\right]$ as the mean first passage time from state $i$ to state $j$ (or mean recurrence time of state $i$ when $i = j$). It is well known that for finite irreducible M. C.'s all the $m_{ij}$ are well defined and finite. Let $M = [m_{ij}]$ be the mean first passage time matrix. Let $\delta_{ij} = 1$, when $i = j$ and 0, when $i \neq j$. Let $M_d = \left[\delta_{ij} m_{ij}\right]$ be the diagonal matrix formed from the diagonal elements of $M$, and $E = [1]$ (i.e. all the elements are unity).

It is well known (Kemeny and Snell (1960)) that, for $1 \leq i, j \leq m$,

$$m_{ij} = 1 + \sum_{k \neq j} p_{ik} m_{kj}. \tag{1.3}$$

In particular, the mean recurrence time of state $j$ is given by

$$m_{jj} = 1/\pi_j. \tag{1.4}$$

From (1.3) and (1.4) it follows that $M$ satisfies the matrix equation

$$(I - P)M = E - PM_d, \text{ with } M_d = \left(\Pi_d\right)^{-1}. \tag{1.5}$$

Generalized matrix inverses (g-inverses) of $I - P$ are typically used to solve systems of linear equations (e.g. (1.2) and (1.5)). Various properties of the M. C., in particular the $\{\pi_j\}$ and the $\{m_{ij}\}$, can be found in terms of g-inverses of $I - P$, either in matrix or elemental form.

While it is possible to use any one-condition g-inverse of $A = I - P$ to solve (1.2) and (1.5), special g-inverses are often used because of their desirable additional properties.

One such g-inverse is $A^\#$, the "*group inverse*" of the matrix $A = I - P$, which is the unique matrix satisfying not only the condition $AA^\# A = A$, but also the additional conditions $A^\# AA^\# = A^\#$ and $AA^\# = A^\# A$. When the M. C. is ergodic, $A^\#$ has the representation $A^\# \equiv [I - P + \Pi]^{-1} - \Pi$, as originally identified by Meyer (1975).

The group inverse $A^\#$, of $A = I - P$, has a number of important properties and leads to modern efficient methods for analysing M. C.'s, (Berman and Plemmons (1979, 1994)). In particular, Meyer (1975), shows that

$$I - AA^{\#} = \begin{cases} \lim_{n\to\infty} \dfrac{I + P + P^2 + ... + P^n}{n}, & \text{for every M. C.,} \\ \lim_{n\to\infty} \sum_{i=0}^{n} \begin{pmatrix} n \\ i \end{pmatrix} k^{n-i}(1-k)^i P^i, & \text{if the M. C. is ergodic and } 0 < k < 1, \\ \lim_{n\to\infty} P^n, & \text{if the M. C, is either regular or absorbing.} \end{cases}$$

We are interested in developing techniques for finding three key properties of discrete time finite M. C.'s: (*i*) the stationary probabilities $\{\pi_j\}$, $(1 \le j \le m)$, (*ii*) the mean first passage times $\{m_{ij}\}$ $(1 \le i, j \le m)$, and (*iii*) the group inverse of $A = I - P$, $A^{\#}$. These properties provide significant information regarding (*i*) the long-term behaviour, (*ii*) the short-term behaviour, and (*iii*) the key attributes of the M. C. The key focus in this paper is on using perturbation techniques in an attempt to develop some useful procedures. Langville and Meyer (2006) pointed out that sequential rank-one updating algorithms typically cost $O(m^3)$ flops and therefore may not be practical in general. While the algorithms may not be computationally efficient it is useful to explore whether we can obtain accurate results. We should point out that determining stationary distributions for M. C.'s using perturbations has been considered by a number of researchers (see, for example, Meyer (1980), Meyer and Shoaf (1980), Funderlic and Plemmons (1986), Hunter (1986), Seneta (1991), Hunter (1991), Langville and Meyer (2006)). However it is the extension to consider the mean first passage times and the group inverse that is a main aim of this paper.

Before exploring this we focus first on generalized matrix inverses, their properties and applications in expressions for the key properties of M. C.'s.

Following a discussion on computational considerations, we describe four different algorithms, all based upon perturbation procedures, where we make row-by-row changes to the transition matrix. We highlight the computational differences between these different algorithms by using a typical five state M. C. to make comparisons.

A sequel to this paper is planned by comparing the perturbation algorithms of this paper, with current techniques for finding the mean first passage times, as well as an alternative accurate computational procedure given in Hunter (2016) based upon an algorithm, due to Kohlas (1986). In this last procedure no subtractions need be carried out.

## 2. Generalized Matrix Inverses

We summarize the definition and classification of generalized matrix inverses. For the context of this paper, we restrict attention to real square matrices of dimension $m$.

**Definition 1**: Let $A$ be an $m \times m$ matrix of real elements. Let $X$ be any $m \times m$ matrix such that $X$ satisfies some of the following conditions:

(Condition 1) $AXA = A$.
(Condition 2) $XAX = X$.
(Condition 3) $(AX)^T = AX$.
(Condition 4) $(XA)^T = XA$.
(Condition 5) $AX = XA$.

If $A^{(i,j,\ldots,l)}$ is any matrix $X$ that satisfies conditions (*i*), (*j*), ..., (*l*) of the above itemised conditions, then $A^{(i,j,\ldots,l)}$ is called an - (*i, j, ..., l*) g-inverse of *A*, under the assumption that condition (1) is always included. Let $A\{i, j, \ldots, l\}$ be the class of all (*i, j, ..., l*) g-inverses of *A*.

$A^{(1)}$, a one-condition g-inverse of *A*, is often written as $A^-$. If *A* is a non-singular matrix (i.e. *det* (*A*) ≠ 0) then $X = A^{-1}$, the inverse of *A*, satisfies all the conditions of the Definition1 and is unique. If *A* is singular, as is the case for $A = I - P$, $A^-$ is not, in general unique. Special cases include $A^{(1,2)}$, a 'reflexive' g-inverse; $A^{(1,3)}$, a 'least squares' g-inverse; $A^{(1,4)}$, a 'minimum norm' g-inverse; $A^{(1,2,3,4)}$, the unique 'Moore-Penrose' g-inverse; and $A^{(1,2,5)}$, the 'group inverse', which exists and is unique if rank($A$) = rank($A^2$).

The following theorem, due to Hunter (1982), gives a procedure for finding all one-condition g-inverses of $I - P$.

**Theorem 1**: *Let P be the transition matrix of a finite irreducible M. C. with m states and stationary probability vector* $\boldsymbol{\pi}^T = (\pi_1, \pi_2, \ldots, \pi_m)$. *Let* $\boldsymbol{e}^T = (1, 1, \ldots, 1)$ *and* $\boldsymbol{t}$ *and* $\boldsymbol{u}$ *be any vectors.*

*(a)* $I - P + \boldsymbol{t}\boldsymbol{u}^T$ *is non-singular if and only if* $\boldsymbol{\pi}^T\boldsymbol{t} \neq 0$ *and* $\boldsymbol{u}^T\boldsymbol{e} \neq 0$.

*(b) If* $\boldsymbol{\pi}^T\boldsymbol{t} \neq 0$ *and* $\boldsymbol{u}^T\boldsymbol{e} \neq 0$ *then* $[I - P + \boldsymbol{t}\boldsymbol{u}^T]^{-1}$ *is a one-condition g-inverse of* $I - P$.

*(c) All one-condition g-inverses of* $I - P$ *can be expressed as*

$A^- = [I - P + \boldsymbol{t}\boldsymbol{u}^T]^{-1} + \boldsymbol{e}\boldsymbol{f}^T + \boldsymbol{g}\boldsymbol{\pi}^T$ *for arbitrary vectors* $\boldsymbol{f}$ *and* $\boldsymbol{g}$.

□

A useful by-product of the proof of the above theorem are the following results:

$$[I - P + \boldsymbol{t}\boldsymbol{u}^T]^{-1}\boldsymbol{t} = \frac{\boldsymbol{e}}{\boldsymbol{u}^T\boldsymbol{e}}. \tag{2.1}$$

$$\boldsymbol{u}^T[I - P + \boldsymbol{t}\boldsymbol{u}^T]^{-1} = \frac{\boldsymbol{\pi}^T}{\boldsymbol{\pi}^T\boldsymbol{t}}. \tag{2.2}$$

A summary of the parametric forms for all generalized inverses of $I - P$, given in Hunter (1990) and extended in Hunter (2014), follows below.

**Theorem 2**: *If G is any g-inverse of* $I - P$, *where P is the transition matrix of a finite irreducible M. C. with stationary probability vector* $\boldsymbol{\pi}^T$, *then G can be uniquely expressed in parametric form as*

$$G \equiv G(\boldsymbol{\alpha}, \boldsymbol{\beta}, \gamma) = [I - P + \boldsymbol{\alpha}\boldsymbol{\beta}^T]^{-1} + \gamma\boldsymbol{e}\boldsymbol{\pi}^T, \tag{2.3}$$

*where* $\boldsymbol{\alpha}$, $\boldsymbol{\beta}$, and $\gamma$ *involve* $2m - 1$ *independent parameters with the property that*

$$\boldsymbol{\pi}^T\boldsymbol{\alpha} = 1 \text{ and } \boldsymbol{\beta}^T\boldsymbol{e} = 1. \tag{2.4}$$

*If* $A_G \equiv I - (I - P)G$ *and* $B_G \equiv I - G(I - P)$ *then*

$$A_G = \boldsymbol{\alpha}\boldsymbol{\pi}^T \text{ and } B_G = \boldsymbol{e}\boldsymbol{\beta}^T, \tag{2.5}$$

*so that* $\boldsymbol{\alpha} = A_G\boldsymbol{e}$ *and* $\boldsymbol{\beta}^T = \boldsymbol{\pi}^T B_G$. (2.6)

*Further, from* (2.1) *and* (2.2), $G\boldsymbol{\alpha} = (\gamma + 1)\boldsymbol{e}$ *and* $\boldsymbol{\beta}^T G = (\gamma + 1)\boldsymbol{\pi}^T$, (2.7)

*so that* $\gamma + 1 = \boldsymbol{\pi}^T G\boldsymbol{\alpha} = \boldsymbol{\beta}^T G\boldsymbol{e} = \boldsymbol{\beta}^T G\boldsymbol{\alpha}$. (2.8)

*Also*

$$G \in A\{1, 2\} \Leftrightarrow \gamma = -1, \quad (2.9)$$

$$G \in A\{1, 3\} \Leftrightarrow \boldsymbol{\alpha} = \boldsymbol{\pi} / \boldsymbol{\pi}^T\boldsymbol{\pi}, \quad (2.10)$$

$$G \in A\{1, 4\} \Leftrightarrow \boldsymbol{\beta} = \boldsymbol{e} / \boldsymbol{e}^T\boldsymbol{e} = \boldsymbol{e} / m, \quad (2.11)$$

$$G \in A\{1, 5a\} \Leftrightarrow \boldsymbol{\alpha} = \boldsymbol{e}, \quad (2.12)$$

$$G \in A\{1, 5b\} \Leftrightarrow \boldsymbol{\beta} = \boldsymbol{\pi}, \quad (2.13)$$

$$G \in A\{1, 5\} \Leftrightarrow \boldsymbol{\alpha} = \boldsymbol{e}, \boldsymbol{\beta} = \boldsymbol{\pi}. \quad (2.14)$$

□

The following results, given in Hunter (2014), provide simple conditions for determining the $A\{1, 5a\}$ and $A\{1, 5b\}$ classes of g-inverses.

**Theorem 3**: *Let* $G = G(\boldsymbol{\alpha}, \boldsymbol{\beta}, \gamma)$ *be any g-inverse of I – P, where P is the transition matrix of a finite irreducible M. C. with stationary probability vector* $\boldsymbol{\pi}^T$.

*(a)* $G \in A\{1, 5a\} \Leftrightarrow G\boldsymbol{e} = g\boldsymbol{e}$ *for some g. Further, if* $G\boldsymbol{e} = g\boldsymbol{e}$ *for some g then* $g = 1 + \gamma$.

*(b)* $G \in A\{1, 5b\} \Leftrightarrow \boldsymbol{\pi}^T G = h\boldsymbol{\pi}^T$ *for some h. Further, if* $\boldsymbol{\pi}^T G = h\boldsymbol{\pi}^T$ *for some h then* $h = 1 + \gamma$.

*(c) If* $G\boldsymbol{e} = g\boldsymbol{e}$ *for some g, and* $\boldsymbol{\pi}^T G = h\boldsymbol{\pi}^T$ *for some h, then* $g = h = \gamma + 1$ *and consequently* $G = G(\boldsymbol{e}, \boldsymbol{\pi}, \gamma)$ *and* $G \in A\{1, 5\}$.

□

Note that the group inverse $A^{\#}$ of $A = I - P$, when $P$ is irreducible, is the unique member of the $A\{1, 2, 5\}$ class of g-inverses with parametric form $G = G(\boldsymbol{e}, \boldsymbol{\pi}, -1)$ so that, from (2.3), $A^{\#} = [I - P + \Pi]^{-1} - \Pi$ where $\Pi = \boldsymbol{e}\boldsymbol{\pi}^T$. This expression for $A^{\#}$ first appeared in Theorem 5.5 of Meyer (1975). $A^{\#}$ has some special properties that can be deduced from Theorems 2 and 3 above.

**Theorem 4**: *If* $A^{\#}$ *is the group inverse of* $A = I - P$, *where P is the transition matrix of a finite irreducible M. C., then the following four properties uniquely determine* $A^{\#}$:

*(i)* $(I - P)A^{\#} = I - \boldsymbol{e}\boldsymbol{\pi}^T$, *(ii)* $A^{\#}(I - P) = I - \boldsymbol{e}\boldsymbol{\pi}^T$, *(iii)* $A^{\#}\boldsymbol{e} = \boldsymbol{0}$, *(iv)* $\boldsymbol{\pi}^T A^{\#} = \boldsymbol{0}^T$.

□

Note that from Theorem 2, the conditions of Theorem 4 imply that $\boldsymbol{\alpha} = \boldsymbol{e}$, $\boldsymbol{\beta}^T = \boldsymbol{\pi}^T$ and $\gamma = -1$, leading to $A^{\#}$ as the group inverse, the only member of $A\{1, 2, 5\}$.

**Theorem 5**: *Let* $G = G(\boldsymbol{\alpha}, \boldsymbol{\beta}, \gamma)$ *be any one-condition g-inverse of I – P.*

(*a*) *Let* $H = G(I - \Pi)$. *Then H is a one-condition g-inverse of I – P with the property that* $H\boldsymbol{e} = \boldsymbol{0}$ and $H = G(\boldsymbol{e}, \boldsymbol{\beta}, -1) = [I - P + \boldsymbol{e}\boldsymbol{\beta}^T]^{-1} - \boldsymbol{e}\boldsymbol{\pi}^T$, *implying that* $H \in A\{1, 2, 5a\}$.

*(b) Let* $K = (I - \Pi)G(I - \Pi) = (I - \Pi)H$. *Then K is an invariant g-inverse and takes the value of the group inverse* $A^{\#}$. *i.e.* $K = G(\boldsymbol{e}, \boldsymbol{\pi}, -1) = [I - P + \boldsymbol{e}\boldsymbol{\pi}^T]^{-1} - \boldsymbol{e}\boldsymbol{\pi}^T$ *implying that* $K \in A\{1, 2, 5\}$.

□

Theorem 5(a) is given by Theorem 7 of Hunter (2014) while Theorem 5(b) appeared in Theorem 6.3 of Hunter (1982) and Corollary 4.6.1 of Hunter (1988).

Thus $H = G(I - \Pi)$ has some of the characteristics ($\boldsymbol{\alpha} = \boldsymbol{e}$ and $\gamma = -1$) of $A^{\#}$ but not all. The additional computation $K = (I - \Pi)H$ in fact characterises $A^{\#}$. The result (b) of Theorem 5 is a useful tool for finding the group inverse when any one-condition g-inverse is available.

One other important related matrix associated with finite irreducible M. C.'s is Kemeny and Snell's *fundamental matrix* $Z = [I - P + \Pi]^{-1}$ which was introduced in Kemeny & Snell (1960). $Z$ was shown in Hunter (1969) to be a one-condition g-inverse of $I - P$ and further, in Hunter (1988), to be a (1, 5) g-inverse with the characterisation $Z = G(\boldsymbol{e}, \boldsymbol{\pi}, 0)$.

One-condition generalized matrix inverses play a major role in solving systems of linear equations often leading to them being called "equation solving" g-inverses. See Theorem 3.1 of Hunter (1982).

## 3. Stationary Distributions using generalized matrix inverses

The papers Hunter (1982) and Hunter (2007) give a variety of expressions for the stationary probability vector $\boldsymbol{\pi}^T$ in terms of different g-inverses of $I - P$. In particular we list the following results that have relevance to the algorithms to be developed in this paper. Let $\boldsymbol{e}_i^T$ be the *i-th* elementary row vector, with 1 in the *i-th* position and 0 elsewhere.

**Theorem 6**:
(*a*) (Hunter, (1982)). *If G is any g-inverse of* $I - P$, $A_G \equiv I - (I - P)G$ *and* $\boldsymbol{v}^T$ *is any vector such that* $\boldsymbol{v}^T A_G \boldsymbol{e} \neq 0$ *then*

$$\boldsymbol{\pi}^T = \frac{\boldsymbol{v}^T A_G}{\boldsymbol{v}^T A_G \boldsymbol{e}}, \tag{3.1}$$

*Furthermore* $A_G \boldsymbol{e} \neq 0$ *for all g-inverses of G, so that it is always possible to find a suitable* $\boldsymbol{v}^T$.

*(b)* (Hunter, (1992)). *If G is a* $(1, 5)$ *g-inverse of* $I - P$,

$$\boldsymbol{\pi}^T = \frac{\boldsymbol{e}^T A_G}{\boldsymbol{e}^T \boldsymbol{e}} \text{ and, for any } i = 1, 2, ..., m, \ \boldsymbol{\pi}^T = \boldsymbol{e}_i^T A_G. \tag{3.2}$$

(c) (Paige, Styan Wachter (1975), Hunter (1982)).
*If* $G = [I - P + \boldsymbol{t}\boldsymbol{u}^T]^{-1}$ *where* $\boldsymbol{u}$ *and* $\boldsymbol{t}$ *are any vectors such that* $\boldsymbol{\pi}^T \boldsymbol{t} \neq 0$ *and* $\boldsymbol{u}^T \boldsymbol{e} \neq 0$, *then*

$$\boldsymbol{\pi}^T = \frac{\boldsymbol{u}^T G}{\boldsymbol{u}^T G \boldsymbol{e}}. \tag{3.3}$$

*In particular*

(*i*) *If* $G = [I - P + \boldsymbol{e}\boldsymbol{u}^T]^{-1}$ *where* $\boldsymbol{u}^T \boldsymbol{e} \neq 0$, *then* $\boldsymbol{\pi}^T = \boldsymbol{u}^T G$. (3.4)

*(ii) If* $G \equiv G_{eb} = [I - P + \boldsymbol{e}\boldsymbol{e}_b^T]^{-1} = [g_{ij}]$ *then* $\pi_j = g_{bj}$. (3.5)

*(iii) If* $G \equiv G_{ee} = [I - P + \boldsymbol{e}\boldsymbol{e}^T]^{-1} = [g_{ij}]$ *then* $\pi_j = \sum_{k=1}^{m} g_{kj} = g_{\bullet j}$. (3.6)

□

If one wishes to find a computationally efficient algorithm for finding $\pi_j$ based upon (3.4), note that we need to solve the equations $\boldsymbol{\pi}^T(I - P + \boldsymbol{e}\boldsymbol{u}^T) = \boldsymbol{u}^T$. Paige, Styan and Wachter (1975) recommended solving this system of linear equations for $\boldsymbol{\pi}$ with $\boldsymbol{u}^T = \boldsymbol{e}_j^T P = \boldsymbol{p}_j^{(r)T}$ , using Gaussian elimination with pivoting.

A numerically stable algorithm for finding the stationary probabilities is the *GTH*/State Reduction algorithm of Grassman, Taksar and Heyman (1985) and Sheskin (1985) that has the advantage that no subtractions are required.

## 4. Mean first passage times using generalized matrix inverses

The solution of equations of the form of (1.5) can be implemented using g-inverses of $I - P$.

**Theorem 7**:
*Let G be any g-inverse of I – P, let* $D = (\Pi_d)^{-1}$ *, and let M be the mean first passage time matrix. Then*
*(a)* $M = [G\Pi - E(G\Pi)_d + I - G + EG_d]D,$ (4.1)
*(b) If* $H \equiv G(I - \Pi)$ *then H is a g-inverse of I – P with* $H\boldsymbol{e} = \boldsymbol{0}$ *and*
$M = [I - H + EH_d]D,$ (4.2)
*(c)* $M = [I - G + EG_d]D,$ (4.3)
*if and only if* $G\boldsymbol{e} = g\boldsymbol{e}$ *for some g, (or equivalently that* $G \in A\{1, 5a\}$*).* □

Equation (4.1) was derived in Hunter (1982). The special case given by (4.2) appears in Hunter (2006) while (4.3) appears in Hunter (2014).

The advantage of the alternative expressions (4.2) and (4.3) is that they lead to simpler elemental forms for the $m_{ij}$, as summarised below.

**Corollary 7.1**:
*Let* $G = [g_{ij}]$ *be any g-inverse of* $I - P$*, and* $M = [m_{ij}]$.
*Let* $g_{i.} = \sum_{j=1}^{m} g_{ij}$ *and* $H = G(I - \Pi) = [h_{ij}]$ *so that*

$$h_{ij} = g_{ij} - g_{i.}\pi_j \text{ for all } i, j \quad (4.4)$$

and

$$m_{ij} = [h_{jj} - h_{ij} + \delta_{ij}]/\pi_j \text{ for all } i, j. \quad (4.5)$$

*In addition, if* $g_{i.} = g$ *, or equivalently that G is a (1, 5a) g-inverse, then*

$$m_{ij} = [g_{jj} - g_{ij} + \delta_{ij}]/\pi_j, \text{ for all } i, j,$$

$$= \begin{cases} [g_{jj} - g_{ij}]/\pi_j, & i \neq j, \\ 1/\pi_j, & i = j. \end{cases} \quad (4.6)$$

□

Special cases of equation (4.3) for $M$ and (4.6) for the elements $m_{ij}$ are $G = Z$, Kemeny and Snell's fundamental matrix, (Kemeny and Snell (1960)) and $G = A^{\#}$, Meyer's group inverse of $I - P$, (Meyer (1975)).

**Corollary 7.2:** (Hunter (2007)).
*If* $G_{eb} = [\, g_{ij}\,] = [I - P + \boldsymbol{e}\boldsymbol{e}_b^T]^{-1}$ *then*

$$\pi_j = g_{bj},\ j = 1,2,\ \ldots, m, \tag{4.7}$$

*with* $m_{ij}$ *given by* (4.6).

□

Thus following one matrix inversion (actually only the *b-th* row, typically the first row, for the stationary distribution), one can find the stationary probabilities and the mean first passage times. We explore this further in some of our results to follow.

We have seen that the mean first passage times $m_{ij}$ can be found using $A^{\#}$. Conversely $A^{\#}$ can in fact be found directly from the $m_{ij}$. The following result appears in Ben Ari & Neumann (2012) and Hunter (2014).

**Theorem 8**: Let $\tau_j = \sum_{k=1}^{m} \pi_k m_{kj} = \sum_{k \neq j} \pi_k m_{kj} + 1,$ *and let* $A^{\#} = [a_{ij}^{\#}]$ *then*

$$a_{ij}^{\#} = \begin{cases} \pi_j(\tau_j - 1), & i = j, \\ \pi_j(\tau_j - 1 - m_{ij}) = a_{jj}^{\#} - \pi_j m_{ij}, & i \neq j. \end{cases}$$

□

Heyman and Reeves (1989) noted that the computation of mean first passage times using the group inverse and the relevant equation (4.3), viz. $M = [I - A^{\#} + EA^{\#}_d]D$ leads to a significant inaccuracy on the more difficult problems in that the computation of $M$ yields three sources of errors:

1. The algorithm for computing $\boldsymbol{\pi}^T$,
2. The computation of the inverse of $I - P + \Pi$, as the matrix may have negative elements that can cause round-off errors in computing the inverse.
3. The matrix evaluation of $M$, as the matrix multiplying $D$ may have negative elements.

Heyman and O'Leary (1995) state that … "it does not make sense to compute …. the group generalized inverse unless the individual elements of those matrices are of interest."

In a recent paper, (Hunter (2016)), an accurate procedure for computing the matrix of mean first passage times is given, based on the procedure of Kohlas (1986). It is shown in that paper that the more general setting of Markov renewal processes leads to a procedure not involving any subtractions. A further paper comparing the results of this paper with other alternative procedures is planned, including this new procedure.

## 5. Perturbed Markov chains

We now explore some relationships between the stationary distributions of unperturbed and perturbed M. C.'s utilizing some special g-inverses of $I - P$. The following results appear in Section 8 of Hunter (2014).

**Theorem 9**: *Let* $P$ *be the transition matrix of a finite irreducible M. C. Let* $\bar{P} = P + \boldsymbol{E}$, *be the transition matrix of the perturbed M. C., with* $\boldsymbol{E} = [\varepsilon_{ij}]$ *as the matrix of perturbations* (*with* $\sum_{j=1}^{m} \varepsilon_{ij} = 0$) *so that* $\boldsymbol{Ee} = \boldsymbol{0}$. *Assume that* $\bar{P}$ *is irreducible. Let* $\boldsymbol{\pi}^T = (\pi_1, \pi_2, \ldots, \pi_m)$ *and* $\bar{\boldsymbol{\pi}}^T = (\bar{\pi}_1, \bar{\pi}_2, \ldots, \bar{\pi}_m)$ *be the stationary probability vectors of the M. C.'s with transition matrices P and* $\bar{P}$, *respectively. Let* $\Pi = \boldsymbol{e\pi}^T$ *and let G be any g-inverse of I – P.*

*(a) If H = G(I – Π), then* $\bar{\boldsymbol{\pi}}^T - \boldsymbol{\pi}^T = \bar{\boldsymbol{\pi}}^T \boldsymbol{E}H.$ (5.1)

*(b) If* $G\boldsymbol{e} = g\boldsymbol{e}$ *for some g, then* $\bar{\boldsymbol{\pi}}^T - \boldsymbol{\pi}^T = \bar{\boldsymbol{\pi}}^T \boldsymbol{E}G.$ (5.2)

*(c) If* $G = [I - P + \boldsymbol{eu}^T]^{-1} + \boldsymbol{ef}^T + \mathbf{g}\boldsymbol{\pi}^T$ *with* $\boldsymbol{u}^T\boldsymbol{e} \neq 0$, $\boldsymbol{f}^T$ *and* $\mathbf{g}$ *arbitrary vectors, then*

$$\bar{\boldsymbol{\pi}}^T - \boldsymbol{\pi}^T = \bar{\boldsymbol{\pi}}^T \boldsymbol{E}[I - P + \boldsymbol{eu}^T]^{-1}. \quad (5.3)$$

*(d) If* $G = A^{\#}$, *the group inverse of I – P then* $\bar{\boldsymbol{\pi}}^T - \boldsymbol{\pi}^T = \bar{\boldsymbol{\pi}}^T \boldsymbol{E}A^{\#}.$ (5.4)

□

If we can establish conditions under which matrices of the form $[I - \boldsymbol{E}H]^{-1}$ exist (or more simply with $H$ taken as $G$) and are of simple form then we can establish useful expressions for $\bar{\boldsymbol{\pi}}^T$ from equations (5.1) to (5.4)). We summarise some special cases where such expressions do in fact hold.

**Theorem 10**: *Under the conditions of Theorem 9,*

(*a*) $I - \boldsymbol{E}A^{\#}$ *is non singular and*

$$\bar{\boldsymbol{\pi}}^T = \boldsymbol{\pi}^T (I - \boldsymbol{E}A^{\#})^{-1}. \quad (5.5)$$

*(b) If* $\mathbf{E} = \boldsymbol{ab}^T$, *let* $\boldsymbol{b}^T H \equiv \boldsymbol{h}^T$ *and assume* $\boldsymbol{h}^T\boldsymbol{a} \neq 1$. *Then* $I - \boldsymbol{E}H = I - \boldsymbol{ah}^T$ *is non-singular and*

$$\bar{\boldsymbol{\pi}}^T = \boldsymbol{\pi}^T [I - \boldsymbol{E}H]^{-1} = \boldsymbol{\pi}^T \left( I + \frac{\boldsymbol{ah}^T}{1 - \boldsymbol{h}^T\boldsymbol{a}} \right). \quad (5.6)$$

*(c) If* $G\boldsymbol{e} = g\boldsymbol{e}$ *for some g,* $\mathbf{E} = \boldsymbol{ab}^T$, *let* $\boldsymbol{b}^T G \equiv \boldsymbol{g}^T$. *Then* $I - \mathbf{E}G = I - \boldsymbol{ag}^T$ *is non-singular and*

$$\bar{\boldsymbol{\pi}}^T = \boldsymbol{\pi}^T [I - \boldsymbol{E}G]^{-1} = \boldsymbol{\pi}^T \left( I + \frac{\boldsymbol{ag}^T}{1 - \boldsymbol{g}^T\boldsymbol{a}} \right). \quad (5.7)$$

***Proof:***

*(a)* The non-singularity of $I - \boldsymbol{E}A^{\#}$ was established in Theorem 3.1 of Meyer (1980) with (5.5) following from (5.4).

(*b*) The Sherman-Morrison (1949) formula states that

$$\left[I - \boldsymbol{ah}^T\right]\left[I + \frac{\boldsymbol{ah}^T}{1 - \boldsymbol{h}^T\boldsymbol{a}}\right] = I = \left[I + \frac{\boldsymbol{ah}^T}{1 - \boldsymbol{h}^T\boldsymbol{a}}\right]\left[I - \boldsymbol{ah}^T\right], \quad (5.8)$$

provided $\boldsymbol{h}^T\boldsymbol{a} \neq 1$. By the "matrix determinant lemma" (see Section 12.1, Harville, (1997)), $\det[I - \boldsymbol{ah}^T] = 1 - \boldsymbol{h}^T\boldsymbol{a} \neq 0$ so that the inverse of $I - \boldsymbol{E}H = I - \boldsymbol{ah}^T$ exists and, from (5.8), $[I - \boldsymbol{ah}^T]^{-1} = I + \dfrac{\boldsymbol{ah}^T}{1 - \boldsymbol{h}^T\boldsymbol{a}}$, leading to (5.6).

(*c*) Equation (5.7) will follow from the arguments that led to (5.6) provided we can establish that $I - \boldsymbol{E}G = I - \boldsymbol{ab}^T G$ is non-singular where $G$ is any (1, 5a) g-inverse of $A = I - P$. The key to establishing (5.7) is Meyer's result regarding the non-singularity of $I - \boldsymbol{E}A^{\#}$. Now $\det(I - \boldsymbol{E}A^{\#}) = \det(I - \boldsymbol{ab}^T A^{\#}) = 1 - \boldsymbol{b}^T A^{\#}\boldsymbol{a} = 1 - \boldsymbol{b}^T[I - P + \boldsymbol{e\pi}^T]^{-1}\boldsymbol{a} \neq 0.$

Now $G$ has the form $G = [I - P + \boldsymbol{e\beta}^T]^{-1} + \gamma\boldsymbol{e\pi}^T$ where $\boldsymbol{\beta}^T\boldsymbol{e} = 1$, (Theorem 2, Hunter (2014).) Further, from equation (3.6) of Theorem 3.3 in Hunter (1988),

$[I - P + \boldsymbol{e\beta}^T]^{-1} = [I - \boldsymbol{e\beta}^T][I - P + \boldsymbol{e\pi}^T]^{-1} + \boldsymbol{e\pi}^T$ so that since $\boldsymbol{b}^T\boldsymbol{e} = 0$ it is easily seen that

$\det(I - \boldsymbol{E}G) = \det(I - \boldsymbol{ab}^T G) = 1 - \boldsymbol{b}^T G\boldsymbol{a} = 1 - \boldsymbol{b}^T[I - P + \boldsymbol{e\pi}^T]^{-1}\boldsymbol{a} \neq 0.$ □

**Theorem 11**: *Under the conditions of Theorem 9, if* $A = I - P$, $\overline{A} = I - \overline{P}$ with $\overline{P} = P + \boldsymbol{E}$ *the group inverse* $\overline{A}^{\#}$ *of* $\overline{A} = I - \overline{P} = I - P - \boldsymbol{E} = A - \boldsymbol{E}$ *is given by*

$$\overline{A}^{\#} = A^{\#}(I - \boldsymbol{E}A^{\#})^{-1} - \Pi(I - \boldsymbol{E}A^{\#})^{-1}A^{\#}(I - \boldsymbol{E}A^{\#})^{-1}. \tag{5.9}$$

***Proof:***

From Theorem 3.1 of Meyer (1980), the group inverse of $\overline{A}^{\#} = (I - \overline{P})^{\#} = (I - P - \boldsymbol{E})^{\#}$ is given by

$$\overline{A}^{\#} = A^{\#} + A^{\#}\boldsymbol{E}A^{\#}(I - \boldsymbol{E}A^{\#})^{-1} - \Pi(I - \boldsymbol{E}A^{\#})^{-1}A^{\#}(I - \boldsymbol{E}A^{\#})^{-1}. \tag{5.10}$$

This expression is also given in Theorem 5.3.30 of Kirkland and Neumann (2013).
Equation (5.9) follows from (5.10) by verifying that $A^{\#} + A^{\#}\boldsymbol{E}A^{\#}(I - \boldsymbol{E}A^{\#})^{-1} = A^{\#}(I - \boldsymbol{E}A^{\#})^{-1}$.

Alternatively, let $G$ be the expression for $\overline{A}^{\#}$ given by (5.9). With $\overline{A} = A - \boldsymbol{E}$ , it can be shown, following some algebraic manipulations, (upon observing that $\boldsymbol{Ee} = \boldsymbol{0}$ , $A\boldsymbol{e} = \boldsymbol{0}, \boldsymbol{E}\Pi = 0$ and $A\Pi = 0$) that $\overline{A}G = I - \boldsymbol{e}\boldsymbol{\pi}^T = G\overline{A}$, from which $\overline{A}G\overline{A} = \overline{A}$ and $G\overline{A}G = G$ showing that $G$ satisfies the three conditions specified by Definition 1 to be the group inverse of $\overline{A}$. □

We now specialize our results to the case when the perturbing matrix has only one non-zero row, the *i-th* row. $\boldsymbol{b}_i^T$.

**Theorem 12**: *Let P and* $\overline{P}$ *be the transition matrices of finite irreducible M. C.'s with* $\overline{P}$ *differing from P only in the i-th row, so that* $\overline{P} = P + \boldsymbol{e}_i\boldsymbol{b}_i^T$ *for some vector* $\boldsymbol{b}_i^T$ *such that* $\boldsymbol{b}_i^T\boldsymbol{e} = 0$. *Let* $\boldsymbol{\pi}^T$ *and* $\overline{\boldsymbol{\pi}}^T$ *be the stationary probability vectors of the respective M. C.'s. Let G be any one-condition g-inverse of I – P. Let* $H = G(I - \boldsymbol{e}\pi^T)$ *and* $A^{\#}$ *be the group inverse of I – P. Then*

*(a)* *(i)* $$\overline{\boldsymbol{\pi}}^T = \boldsymbol{\pi}^T\left[I + \frac{1}{h_i}\boldsymbol{e}_i\boldsymbol{b}_i^T H\right] \text{ where } h_i = 1 - \boldsymbol{b}_i^T H\boldsymbol{e}_i \neq 0, \tag{5.11}$$

*(ii)* $$\overline{\boldsymbol{\pi}}^T = \boldsymbol{\pi}^T\left[I + \frac{1}{g_i}\boldsymbol{e}_i\boldsymbol{b}_i^T G\right] \text{ where } g_i = 1 - \boldsymbol{b}_i^T G\boldsymbol{e}_i \text{ when } G\boldsymbol{e} = g\boldsymbol{e} \text{ for some } g, \tag{5.12}$$

*(iii)* $$\overline{\boldsymbol{\pi}}^T = \boldsymbol{\pi}^T\left[I + \frac{1}{a_i}\boldsymbol{e}_i\boldsymbol{b}_i^T A^{\#}\right] \text{where } a_i = 1 - \boldsymbol{b}_i^T A^{\#}\boldsymbol{e}_i. \tag{5.13}$$

*(b)* *If* $A^{\#}$ *and* $\overline{A}^{\#}$ *are the group inverses of* $A = I - P$ *and* $\overline{A} = I - \overline{P}$ *then*

$$\overline{A}^{\#} = A^{\#} + \frac{1}{a_i}A^{\#}\boldsymbol{e}_i\boldsymbol{b}_i^T A^{\#} - \frac{\pi_i}{a_i}\boldsymbol{e}\boldsymbol{b}_i^T\left(A^{\#} + \frac{\boldsymbol{b}_i^T(A^{\#})^2\boldsymbol{e}_i}{a_i}I\right)A^{\#}. \tag{5.14}$$

***Proof:***

*(a)* Set $\boldsymbol{E} = \boldsymbol{e}_i\boldsymbol{b}_i^T$ so that $\boldsymbol{a} = \boldsymbol{e}_i$. (*i*) Equation (5.11) follows from (5.6) with $\boldsymbol{h}^T = \boldsymbol{b}_i^T H$. (*ii*) Equation (5.12) follows from (5.7) with $\boldsymbol{g}^T = \boldsymbol{b}_i^T G$. (*iii*) Equation (5.13) follows from (5.5) and (5.8) with $\boldsymbol{g}^T = \boldsymbol{b}_i^T A^{\#}$.

*(b)* From (5.9), $\overline{A}^{\#} = A^{\#}(I - \boldsymbol{e}_i\boldsymbol{b}_i^T A^{\#})^{-1} - \boldsymbol{e}\boldsymbol{\pi}^T(I - \boldsymbol{e}_i\boldsymbol{b}_i^T A^{\#})^{-1}A^{\#}(I - \boldsymbol{e}_i\boldsymbol{b}_i^T A^{\#})^{-1}$.

Now $A^{\#}(I - \boldsymbol{e}_i\boldsymbol{b}_i^T A^{\#})^{-1} = A^{\#} + \frac{1}{a_i} A^{\#}\boldsymbol{e}_i\boldsymbol{b}_i^T A^{\#}$ where $a_i = 1 - \boldsymbol{b}_i^T A^{\#}\boldsymbol{e}_i$ and

$$(I - \boldsymbol{e}_i\boldsymbol{b}_i^T A^{\#})^{-1} A^{\#}(I - \boldsymbol{e}_i\boldsymbol{b}_i^T A^{\#})^{-1} = \left[ I + \frac{1}{a_i}\boldsymbol{e}_i\boldsymbol{b}_i^T A^{\#} \right] A^{\#} \left[ I + \frac{1}{a_i}\boldsymbol{e}_i\boldsymbol{b}_i^T A^{\#} \right]$$

implying after simplification and expansion

$$(I - \boldsymbol{e}_i\boldsymbol{b}_i^T A^{\#})^{-1} A^{\#}(I - \boldsymbol{e}_i\boldsymbol{b}_i^T A^{\#})^{-1} = A^{\#} + \frac{A^{\#}\boldsymbol{e}_i\boldsymbol{b}_i^T A^{\#}}{a_i} + \frac{\boldsymbol{e}_i\boldsymbol{b}_i^T (A^{\#})^2}{a_i} + \frac{(\boldsymbol{b}_i^T (A^{\#})^2\boldsymbol{e}_i)\boldsymbol{e}_i\boldsymbol{b}_i^T A^{\#}}{(a_i)^2}.$$

Now $\boldsymbol{\pi}^T A^{\#} = \boldsymbol{0}^T$ so that

$$\overline{A}^{\#} = A^{\#} + \frac{1}{a_i} A^{\#}\boldsymbol{e}_i\boldsymbol{b}_i^T A^{\#} - \boldsymbol{e}\boldsymbol{\pi}^T \left( \frac{\boldsymbol{e}_i\boldsymbol{b}_i^T (A^{\#})^2}{a_i} + \frac{(\boldsymbol{b}_i^T (A^{\#})^2\boldsymbol{e}_i)\boldsymbol{e}_i\boldsymbol{b}_i^T A^{\#}}{(a_i)^2} \right)$$

$$= A^{\#} + \frac{1}{a_i} A^{\#}\boldsymbol{e}_i\boldsymbol{b}_i^T A^{\#} - \boldsymbol{e} \left( \frac{\pi_i\boldsymbol{b}_i^T (A^{\#})^2}{a_{\boldsymbol{i}}} + \frac{\pi_i(\boldsymbol{b}_i^T (A^{\#})^2\boldsymbol{e}_i)\boldsymbol{b}_i^T A^{\#}}{(a_i)^2} \right)$$

leading to (5.14).

□

In the algorithms to follow in the next section we use a very simple procedure. We start with a simple transition matrix $P_0$ with known or easily computed stationary probability vector $\boldsymbol{\pi}_0^T$, mean first passage time matrix $M_0$ and group inverse $A_0^{\#}$ or a simple g-inverse $G_0$. We then sequentially change the transition matrix $P_0$ by replacing the *i-th* row of $P_0$ with the *i-th* row of $P$ (i.e. $\boldsymbol{p}_i^T = \boldsymbol{e}_i^T P$) ($i = 1, 2, \ldots, m$) to obtain $P_i$ ending up with $P_m = P$.

Thus let $P_0 = \sum_{i=1}^{m} \boldsymbol{e}_i\boldsymbol{p}_{(0)i}^T$ so that if $P = \sum_{i=1}^{m} \boldsymbol{e}_i\boldsymbol{p}_i^T$ then $P_i = P_{i-1} + \boldsymbol{e}_i\boldsymbol{b}_i^T$ with $\boldsymbol{b}_i^T = \boldsymbol{p}_i^T - \boldsymbol{p}_{(0)i}^T$, for $i = 1, 2, \ldots, m$. Thus we update $\boldsymbol{\pi}_{i-1}^T$, $M_{i-1}$ and $A_{i-1}^{\#}$(or $G_{i-1}$) to $\boldsymbol{\pi}_i^T$, $M_i$ and $A_i^{\#}$ (or $G_i$) finishing with $\boldsymbol{\pi}_m^T = \boldsymbol{\pi}^T$, $M_m = M$ and $A_m^{\#} = A^{\#}$.

We need to start with an irreducible transition matrix $P_0$ and ensure that each successive transition matrix $P_i$ is also irreducible. The simplest structure is to take $P_0 = \frac{1}{m}\boldsymbol{e}\boldsymbol{e}^T$, implying $\boldsymbol{p}_{(0)i}^T = \frac{\boldsymbol{e}^T}{m}$. This leads to $M_0 = m\boldsymbol{e}\boldsymbol{e}^T$ and $A_0^{\#} = I - \frac{1}{m}\boldsymbol{e}\boldsymbol{e}^T$.

Let $A_{i-1} = I - P_{i-1}$ with the *i-th* and subsequent rows of $P_{i-1}$ all $\boldsymbol{e}^T/m$, so that for the update $\overline{A}$ is taken as $A_i = I - P_i$ with the *i-th* row of $P_i$ taken as the prescribed vector $\boldsymbol{p}_i^T$. This is equivalent to taking $P_i = P_{i-1} + \boldsymbol{e}_i\boldsymbol{b}_i^T$ with $\boldsymbol{b}_i^T = \boldsymbol{p}_i^T - \boldsymbol{e}^T/m$, $i = 1, 2, \ldots, m$.

## 6. **The algorithms**

We consider a variety of techniques.

1. Extend the procedure of Hunter (1991), updating one-condition generalized inverses to find successive stationary probability vectors, to compute the group inverse and mean first passage time matrix.

2. Successive direct perturbation updates of the group inverses of the perturbed matrices, leading to an expression for the stationary distribution and the group inverse (and hence the mean first passage times).
3. Consider an extension to the second procedure through updating using matrix procedures that yield, in tandem, the stationary probability vectors and the group inverses.
4. Three interrelated algorithms, each with different starting conditions, based on updating simple generalised inverses of $I - P_0$ that lead to simple computations for the stationary probabilities, the group inverse and the mean first passage time matrix,

### 6.1 Procedure using successive updating of general g-inverses of *I* – *P*.

This procedure is based upon Hunter (1991). Let $P_0 = \boldsymbol{e}\boldsymbol{e}^T/m$. For $i$ = 1, 2, …, $m$, let $P_i = P_{i-1} + \boldsymbol{e}_i\boldsymbol{b}_i^T$ with $\boldsymbol{b}_i^T = \boldsymbol{p}_i^T - \boldsymbol{e}^T/m$, and let $G_i = [I - P_i + \boldsymbol{t}_i\boldsymbol{u}_i^T]^{-1}$. We update the g-inverse $G_{i-1}$ to $G_i$ successively as follows. Take $\boldsymbol{t}_0 = \boldsymbol{e}$ and $\boldsymbol{u}_0^T = \boldsymbol{e}^T/m$ then $G_0 = [I - P_0 + \boldsymbol{t}_0\boldsymbol{u}_0^T]^{-1} = I$.

First note that $\boldsymbol{u}_0^T\boldsymbol{e} \neq 0$, $\boldsymbol{\pi}_0^T\ \boldsymbol{t}_0 \neq 0$, $\boldsymbol{\pi}_0^T = \dfrac{\boldsymbol{u}_0^T G_0}{\boldsymbol{u}_0^T G_0\boldsymbol{e}} = \boldsymbol{e}^T/m$.

For $i$ = 1, 2, …, $m$, let $\boldsymbol{t}_i = \boldsymbol{e}_i$ and $\boldsymbol{u}_i^T = \boldsymbol{u}_{i-1}^T + \boldsymbol{b}_i^T = \boldsymbol{u}_{i-1}^T + \boldsymbol{p}_i^T - \boldsymbol{e}^T/m$, then

$G_i = [I - P_i + \boldsymbol{t}_i\boldsymbol{u}_i]^{-1} = G_{i-1}[I + (\boldsymbol{e}_{i-1} - \boldsymbol{e}_i)(\boldsymbol{\pi}_{i-1}^T/\boldsymbol{\pi}_{i-1}^T\boldsymbol{e}_i)]$ leading to $\boldsymbol{\pi}_i^T = \dfrac{\boldsymbol{u}_i^T G_i}{\boldsymbol{u}_i^T G_i\boldsymbol{e}}$.

In Hunter (1991) it is shown that $G_i = G_{i-1} + F_{i-1}$ where all the elements in $F_{i-1}$ in rows numbered $i$+1, …, $m$ are all zero. The basic algorithm is as follows.

**Algorithm 1**

*(i)* Let $G_0 = I$, $\boldsymbol{u}_0^T = \boldsymbol{e}^T/m$.

*(ii)* For $i = 1, 2, ..., m$, let $\boldsymbol{p}_i^T = \boldsymbol{e}_i^T P$, $\boldsymbol{u}_i^T = \boldsymbol{u}_{i-1}^T + \boldsymbol{p}_i^T - \boldsymbol{e}^T/m$,

$G_i = G_{i-1} + G_{i-1}(\boldsymbol{e}_{i-1} - \boldsymbol{e}_i)(\boldsymbol{u}_{i-1}^T G_{i-1}/\boldsymbol{u}_{i-1}^T G_{i-1}\boldsymbol{e}_i)$.

*(iii)* At $i = m$, let $G_m = G$ and $\boldsymbol{\pi}^T = \boldsymbol{\pi}_m^T = \dfrac{\boldsymbol{u}_m^T G_m}{\boldsymbol{u}_m^T G_m\boldsymbol{e}}$.

*(iv)* Compute $H = G(I - \boldsymbol{e}\boldsymbol{\pi}^T)$.

*(v)* Compute $A^{\#} = (I - \boldsymbol{e}\boldsymbol{\pi}^T)H$.

*(vi)* Compute $M = [I - H + E(diag(H))]D$ where $E = [1]$ *and* $D = inv[diag(\boldsymbol{e}\boldsymbol{\pi}^T)]$.

Some other simplifications are also possible. For example start with

*(i)* $G_1 = I + (\boldsymbol{e} - \boldsymbol{e}_1)\boldsymbol{e}^T$ and let $\boldsymbol{\alpha}_1^T = \boldsymbol{p}_1^T G_1$.

*(ii)* For $i = 1, ..., m-1$, compute

*(a)* $\boldsymbol{v}_i^T = \dfrac{\boldsymbol{\alpha}_i^T}{\boldsymbol{\alpha}_i^T\boldsymbol{e}_{i+1}}$,

*(b)* Compute the first $i$ rows of $B_i = G_i(\boldsymbol{e}_i - \boldsymbol{e}_{i+1})\boldsymbol{v}_i^T$, with the other entries all 0,

*(c)* Set $G_{i+1} = G_i + B_i$,

(d) Compute $\boldsymbol{\alpha}_{i+1}^T = \boldsymbol{v}_i^T - \boldsymbol{e}^T + \boldsymbol{p}_{i+1}^T G_{i+1}$.

(iii) Compute at the final iteration $\boldsymbol{\pi}^T = \dfrac{\boldsymbol{\alpha}_m^T}{\boldsymbol{\alpha}_m^T \boldsymbol{e}}$.

(iv), (v) and (vi) as above.

The justification for the above modifications are discussed in Hunter (1991).

Note that Seneta (1991) also proposed a similar procedure but updating Kemeny and Snell's fundamental matrix $Z$ by taking $P_0 = ee^T P/m$ with $\pi_0^T = e^T P/m$ and $Z_0 = I$. No numerical comparisons were given.

**6.2 Procedure based on row perturbations of the Group Inverse**

Let us explore (5.14) in this recursive setting, by starting with $A^\#$ as $A_{i-1}^\#$ and $\bar{A}^\#$ as $A_i^\#$,

$$A_i^\# = A_{i-1}^\# + \frac{1}{1-\boldsymbol{b}_i^T A_{i-1}^\# \boldsymbol{e}_i} A_{i-1}^\# \boldsymbol{e}_i \boldsymbol{b}_i^T A_{i-1}^\# + \boldsymbol{e}\boldsymbol{y}_i^T, \qquad (6.1)$$

where $\boldsymbol{y}_i^T = -\left(\dfrac{\pi_i^{(i-1)}}{1-\boldsymbol{b}_i^T A_{i-1}^\# \boldsymbol{e}_i}\right)\boldsymbol{b}_i^T\left(A_{i-1}^\# + \dfrac{\boldsymbol{b}_i^T (A_{i-1}^\#)^2 \boldsymbol{e}_i}{1-\boldsymbol{b}_i^T A_{i-1}^\# \boldsymbol{e}_i} I\right)A_{i-1}^\#$.

Note that since $A_{i-1}{}^\# \boldsymbol{e} = \boldsymbol{0}$ it follows that $\boldsymbol{y}_i^T \boldsymbol{e} = 0$. (See Kirkland and Neumann (2013)).

The full computation of (6.1) requires expressions for the constants $1-\boldsymbol{b}_i^T A_{i-1}^\# \boldsymbol{e}$, $\boldsymbol{b}_i^T (A_{i-1}^\#)^2 \boldsymbol{e}_i$ as well as $\pi_i^{(i-1)}$. Note however from (6.1) that if we express the *i-th* group inverse as $A_i^\# = R_i + \boldsymbol{e}\boldsymbol{y}_i^T$, starting with $A_0^\# = R_0$ and $\boldsymbol{y}_0^T = \boldsymbol{0}^T$, then the recursion (6.1) with any terms of the form $\boldsymbol{e}\boldsymbol{g}_i^T$ omitted, can be expressed, using the observation that $\boldsymbol{b}_i^T \boldsymbol{e} = 0$, as

$$R_i = R_{i-1} + \frac{1}{1-\boldsymbol{b}_i^T R_{i-1} \boldsymbol{e}_i} R_{i-1} \boldsymbol{e}_i \boldsymbol{b}_i^T R_{i-1},$$

with $A_i^\#$ found as $R_i + \boldsymbol{e}\boldsymbol{y}_i^T$ with $\boldsymbol{y}_i^T \boldsymbol{e} = 0$. The determination of the $\boldsymbol{y}_i^T$, in particular when $i = m$, when $A^\# = A_m^\#$, can be determined by requiring $R_i \boldsymbol{e} = \boldsymbol{0}$ and the properties of the group inverse. Since $(I-P)A^\# = I - \boldsymbol{e}\boldsymbol{\pi}^T$ we have that $\boldsymbol{\pi}^T = \boldsymbol{e}_1^T - \boldsymbol{e}_1^T (I-P)R$. Further $\boldsymbol{\pi}^T A = \boldsymbol{0}^T$ implies that $\boldsymbol{y}_m^T = -\boldsymbol{\pi}^T R$ so that $A^\# = (I - \boldsymbol{e}\boldsymbol{\pi}^T)R$.

The procedure is outlined below.

**Algorithm 2**

Start with $P$.

(i) Set $R_0 = I - \boldsymbol{e}\boldsymbol{e}^T/m$.

(ii) For $i = 1, 2, \ldots, m$, let $\boldsymbol{p}_i^T = \boldsymbol{e}_i^T P$, $\boldsymbol{b}_i^T = \boldsymbol{p}_i^T - \boldsymbol{e}^T/m$,

$$R_i = R_{i-1} + \frac{1}{1-\boldsymbol{b}_i^T R_{i-1} \boldsymbol{e}_i} R_{i-1} \boldsymbol{e}_i \boldsymbol{b}_i^T R_{i-1}.$$

*(iii)* Compute $\boldsymbol{\pi}^T = \boldsymbol{e}_1^T - \boldsymbol{e}_1^T(I-P)R_m$.

*(iv)* Compute $A^{\#} = (I - \boldsymbol{e}\boldsymbol{\pi}^T)R_m$.

*(v)* Compute $M = [I - R_m + E(diag(R_m))]D$, where $E$ = [1] and $D = inv[diag(\boldsymbol{e}\boldsymbol{\pi}^T)]$.

Some simplifications to this algorithm are possible, as not all the calculations are required. In (*ii*) note that $R_i = R_{i-1}(I + C_i)$, where $C_i = \frac{1}{k_i}\boldsymbol{e}_i\boldsymbol{b}_i^T R_{i-1}$ and $k_i = 1 - \boldsymbol{b}_i^T R_{i-1}\boldsymbol{e}_i$ so that $C_i$ has all terms zero except in the *ith* row. So that in the *ith* recursion the only terms that are updated are in the first *i* rows with the rows numbered *i*+1. *i*+2, …, *m* remaining unchanged.

### 6.3 Procedure based on updating the group inverse by matrix operations

Rather than focus directly on the expression of the group inverse, observe that, from (5.5) under the perturbation $\boldsymbol{E}$, $\bar{\boldsymbol{\pi}}^T = \boldsymbol{\pi}^T(I - \boldsymbol{E}A^{\#})^{-1}$.

Thus, if $\Pi = \boldsymbol{e}\boldsymbol{\pi}^T$ and $\overline{\Pi} = \boldsymbol{e}\overline{\boldsymbol{\pi}}^T$ then $\overline{\Pi} = \Pi(I - \boldsymbol{E}A^{\#})^{-1}$. (6.2)

Now under the perturbation $\boldsymbol{E} = \boldsymbol{e}_i\boldsymbol{b}_i^T$ to the *i*-th row with $\boldsymbol{b}_i^T\boldsymbol{e} = 0$, yields, as in (5.13),

$$(I - \boldsymbol{E}A^{\#})^{-1} = I + \frac{1}{1 - \boldsymbol{b}_i^T A^{\#}\boldsymbol{e}_i}\boldsymbol{e}_i\boldsymbol{b}_i^T A^{\#} \text{ so that } \overline{\Pi} = \Pi\left[I + \frac{1}{1 - \boldsymbol{b}_i^T A^{\#}\mathbf{e}_i}\boldsymbol{e}_i\boldsymbol{b}^T A^{\#}\right]$$

and, from (5.9) and (6.2), $\overline{A}^{\#} = (I - \overline{\Pi})A^{\#}(I - \boldsymbol{E}A^{\#})^{-1} = (I - \overline{\Pi})A^{\#}\left(I + \frac{1}{1 - \boldsymbol{b}_i^T A^{\#}\boldsymbol{e}_i}\boldsymbol{e}_i\boldsymbol{b}_i^T A^{\#}\right)$.

In the context of successive updating of the group inverse on a row by row basis we have the following procedure.

**Algorithm 3**

*(i)* Let $P_0 = \boldsymbol{e}\boldsymbol{e}^T/m$, implying $\Pi_0 = \boldsymbol{e}\boldsymbol{e}^T/m$, $A_0^{\#} = I - \boldsymbol{e}\boldsymbol{e}^T/m$.

*(ii)* For $i = 1, 2, ..., m$, let $\boldsymbol{p}_i^T = \boldsymbol{e}_i^T P$, $\boldsymbol{b}_i^T = \boldsymbol{p}_i^T - \boldsymbol{e}^T/m$,

$$S_i = I + \frac{1}{1 - \boldsymbol{b}_i^T A_{i-1}^{\#}\boldsymbol{e}_i}\boldsymbol{e}_i\boldsymbol{b}_i^T A_{i-1}^{\#},\quad \Pi_i = \Pi_{i-1}S_i,\quad \mathrm{A}_i^{\#} = (I - \Pi_i)A^{\#}_{i-1}S_i.$$

*(iii)* At $i = m$, let $\mathrm{S} = S_m$ then $\Pi = \Pi_{m-1}S$, $A^{\#} = (I - \Pi)A_{m-1}^{\#}S$.

*(iv)* Compute $M = [I - A^{\#} + EA_d^{\#}]D$, where $E$ = [1] and $D = (\Pi_d)^{-1}$.

### 6.4 Procedures based on updating simple g-inverses of *I* – *P*.

We have seen earlier (Theorem 7(c)) that if we choose a g-inverse *G* of *I* – *P* with the property that $G\boldsymbol{e} = g\boldsymbol{e}$, (*i.e.* $G \in A\{1, 5a\}$) then we have a simple form of the mean first passage time matrix *M* given by eqn. (4.3). Further, it is easy to find an expression for the group inverse of $I - P$ as $A^{\#} = (I - \boldsymbol{e\pi}^T)G$. If we take $G$ of the form $G = [I - P + \boldsymbol{e\beta}^T]^{-1}$, then its computation would not require any prior knowledge of the stationary probability vector $\boldsymbol{\pi}^T$.

In Hunter (2007) we explored the properties of some generalized inverses of this form. We consider three different algorithms using the special forms, $G_e \equiv \left[I - P + \frac{\boldsymbol{ee}^T}{m}\right]^{-1}$, $G_{e1} \equiv [I - P + \boldsymbol{ee}_1^T]^{-1}$ and $G_{ee} \equiv [I - P + \boldsymbol{ee}^T]^{-1}$, utilizing (3.4), (3.5) and (3.6), respectively.

The starting conditions for each algorithm, followed by similar recursions, but with different expressions for the stationary probability vector $\boldsymbol{\pi}^T$ lead to identical calculation procedures for the group inverse and the mean first passage times.

We explore the recursions to determine $G = [I - P + \boldsymbol{e\beta}^T]^{-1}$. In each case we start with $K_0 = [I - P_0 + \boldsymbol{e\beta}^T]^{-1} = \left[I - \frac{\boldsymbol{ee}^T}{m} + \boldsymbol{e\beta}^T\right]^{-1} = \left[I + \boldsymbol{eh}^T\right]^{-1}$ where $\boldsymbol{h}^T = \boldsymbol{\beta}^T - \frac{\boldsymbol{e}^T}{m}$ implying, from the proof of Theorem 10(b), that $K_0 = I - \frac{\boldsymbol{eh}^T}{1 + \boldsymbol{h}^T\boldsymbol{e}}$.

The recursion is to take $K_{i-1} = [I - P_{i-1} + \boldsymbol{e\beta}^T]^{-1}$ to $K_i = [I - P_i + \boldsymbol{e\beta}^T]^{-1}$ where $P_i = P_{i-1} + \boldsymbol{e}_i\boldsymbol{b}_i^T$ and $\boldsymbol{b}_i^T = \boldsymbol{p}_i^T - \boldsymbol{e}^T/m$.

Now $K_i = [I - P_i + \boldsymbol{e\beta}^T]^{-1} = [I - P_{i-1} + \boldsymbol{e\beta}^T - \boldsymbol{e}_i\boldsymbol{b}_i^T]^{-1} = [(K_{i-1})^{-1} - \boldsymbol{e}_i\boldsymbol{b}_i^T]^{-1}$.

Using the well known Sherman-Morrison (1949) formula: If $A$ is invertible and $(A + \boldsymbol{uv}^T)^{-1} = A^{-1} - \frac{1}{1 + \boldsymbol{v}^T A^{-1}\boldsymbol{u}} A^{-1}\boldsymbol{uv}^T A^{-1}$, the above expression leads to the recursion:

for $i = 1, 2, \ldots, m$, $K_i = K_{i-1} + \frac{1}{1 - \boldsymbol{b}_i^T K_{i-1}\boldsymbol{e}_i} K_{i-1}\boldsymbol{e}_i\boldsymbol{b}_i^T K_{i-1}$, (6.3)

with $K_0 = [I - P_0 + \boldsymbol{e\beta}^T]^{-1}$ and $K_m = [I - P_m + \boldsymbol{e\beta}^T]^{-1} = [I - P + \boldsymbol{e\beta}^T]^{-1}$.

For $G_e = \left[I - P + \frac{\boldsymbol{ee}^T}{m}\right]^{-1} = K_m$, $\boldsymbol{\beta}^T = \frac{\boldsymbol{e}^T}{m}$, $K_0 = I$ and $\boldsymbol{\pi}^T = \frac{1}{m}\boldsymbol{e}^T K_m$.

For $G_{e1} = [I - P + \boldsymbol{ee}_1^T]^{-1} = K_m$, $\boldsymbol{\beta}^T = \boldsymbol{e}_1{}^T$, $K_0 = I + \boldsymbol{e}\left(\frac{\boldsymbol{e}^T}{m} - \boldsymbol{e}_1^T\right)$ and $\pi^T = \boldsymbol{e}_1^T K_m$.

For $G_{ee} = [I - P + \boldsymbol{ee}^T]^{-1} = K_m$, $\boldsymbol{\beta}^T = \boldsymbol{e}^T$, $K_0 = I - \left(\frac{m-1}{m^2}\right)\boldsymbol{ee}^T$ and $\boldsymbol{\pi}^T = \boldsymbol{e}^T K_m$.

The expressions for $\boldsymbol{\pi}^T$ follow from (2.2) since $\boldsymbol{\pi}^T = \boldsymbol{\beta}^T[I - P + \boldsymbol{e\beta}^T]^{-1}$.

This leads to three further algorithms, all variants of the generic recursion given by (6.3).

**Algorithm 4**: (and its variants $4A$, $4B$ and $4C$)

($i$) Start with $K_0$ : (For $AL4A$ let $K_0 = I$. For $AL4B$ let $K_0 = I + \boldsymbol{e}\left(\frac{\boldsymbol{e}^T}{m} - \boldsymbol{e}_1^T\right)$.

For $AL4C$ let $K_0 = I - \left(\frac{m-1}{m^2}\right)\boldsymbol{e}\boldsymbol{e}^T$.)

($ii$) For $i = 1, 2, ..., m$, let $\boldsymbol{p}_i^T = \boldsymbol{e}_i^T P$, $\boldsymbol{b}_i^T = \boldsymbol{p}_i^T - \boldsymbol{e}^T/m$,

$$K_i = K_{i-1}\left(I + C_i\right), \textit{where } k_i = 1 - \boldsymbol{b}_i^T K_{i-1}\boldsymbol{e}_i \textit{ and } C_i = \frac{1}{k_i}\boldsymbol{e}_i\boldsymbol{b}_i^T K_{i-1}.$$

($iii$) At $i = m$, let $K = K_m$ and then compute $\boldsymbol{\pi}^T$ :

(For $AL4A$ let $\boldsymbol{\pi}^T = \frac{1}{m}\boldsymbol{e}^T K$. For $AL4B$ let $\pi^T = \boldsymbol{e}_1^T K$. For $AL4C$ let $\boldsymbol{\pi}^T = \boldsymbol{e}^T K$.)

($iv$) Compute $A^{\#} = (I - \boldsymbol{e}\boldsymbol{\pi}^T)K$.

($v$) Compute $M = [I - K + EK_d]D$, where $E = [1]$ and $D = (\Pi_d)^{-1}$.

## 7. Numerical results

We conclude this paper with a comparison of all the algorithms by coding each algorithm using MatLab and exploring numerical computations using a test example (which has previously been considered in the literature). MatLab was run in both single precision and double precision. As had been done by others, a comparison of the single precision and double precision results has been used to compare the accuracy of the different algorithms.

One of the difficulties in making comparisons as to which algorithm is preferable is that there is no bench mark of accurate results related to test problems in respect to the mean first passage times and the group inverse.

### 7.1 Stationary distributions

Using the four algorithms, as listed in section 6, the stationary distributions were computed in single precision giving $\{\pi_i(S)\}$ and double precision $\{\pi_i(D)\}$.

In order to compare the procedures against a benchmark procedure we used the *GTH*/State Reduction algorithm of Grassman, Taksar and Heyman (1985) and Sheskin (1985). This was carried out in single precision and double precision and listed as $\{\pi_i(GTHS)\}$ and $\{\pi_i(GTHD)\}$ respectively. The double precision figures of the *GTH* algorithm were taken as the most accurate that we could obtain. Note that the stationary probability vector is just the left eigenvector of $P$ corresponding to the dominant eigenvalue 1. One could alternatively have used MatLab's *eigs* package as an alternative benchmark. However we elected to use the *GTH* algorithm due to its numerical stability with no subtractions being used in the calculations.

We also calculate the *average number of accurate decimal places* for both the single precision and double precision results for each algorithm by comparing the actual computed results against an appropriately rounded version of the *GTH* computed distribution.

Suppose $\{\pi_i(A)\}$ is the calculation for Procedure $A$ and $\{\pi_i(B)\}$ is the calculation for Procedure $B$. The *minimum error, the maximum error and the relative error between A and B* are, respectively, $MINE(A, B) = \min_{1\le i\le m}|\pi_i(A) - \pi_i(B)|$, $MAXE(A,B) = \max_{1\le i\le m}|\pi_i(A) - \pi_i(B)|$, and $RELE(A,B) = \sum_{i=1}^{m}|\pi_i(A) - \pi_i(B)|$.

Under, single and double precision, for any particular algorithm, the *minimum residual error,* the *maximum residual error* and the *relative errors* are, respectively, $MINRE(\cdot) = \min_{1\le j\le m}\left|\pi_j(\cdot) - \sum_{i=1}^{m}\pi_i(\cdot)p_{ij}\right|$, $MAXRE(\cdot) = \max_{1\le j\le m}\left|\pi_j(\cdot) - \sum_{i=1}^{m}\pi_i(\cdot)p_{ij}\right|$, and $RELE(\cdot)$ = $\sum_{j=1}^{m}\left|\pi_j(\cdot) - \sum_{i=1}^{m}\pi_i(\cdot)p_{ij}\right|$.

We provide a table of comparisons comparing the single and double precision results together with a comparison against the *GTH* algorithm.

## 7.2 Mean first passage times

A measure for the accuracy of the mean first passage times was carried out by calculating the $m_{ij}$ using the algorithms in single and double precision to compute and compare the matrices $M(S) = [m_{ij}(S)]$ and $M(D) = [m_{ij}(D)]$.

We could not find any published results for accurate values of the $m_{ij}$ against any specific test problem in the literature.

We consider the *minimum*, *maximum and overall residual errors* for each algorithm (based on the formal calculation for the $m_{ij}$ given by (1.3)), under both single and double precision. i.e $MINRESM(\cdot) = \min_{1\le i\le m,\ 1\le j\le m}\left|m_{ij}(\cdot) - \sum_{k\ne j}p_{ik}m_{kj}(\cdot) - 1\right|$, $MAXRESM(\cdot) = \max_{1\le i\le m,\ 1\le j\le m}\left|m_{ij}(\cdot) - \sum_{k\ne j}p_{ik}m_{kj}(\cdot) - 1\right|$, and $RESM(\cdot) = \sum_{i=1}^{m}\sum_{j=1}^{m}\left|m_{ij}(\cdot) - \sum_{k\ne j}p_{ik}m_{kj}(\cdot) - 1\right|$.

The accuracy of each algorithm was evaluated in terms of the *minimum error, the maximum error and the relative errors between the double and single precision computations as* $MINEM(S, D) = \min_{1\le i\le m,\ 1\le j\le m}|m_{ij}(S) - m_{ij}(D)|$, $MAXEM(S, D) = \max_{1\le i\le m,\ 1\le j\le m}|m_{ij}(S) - m_{ij}(D)|$ and $REM(S,D) = \sum_{i=1}^{m}\sum_{j=1}^{m}|m_{ij}(S) - m_{ij}(D)|$.

If one regards the double precision result as the "true" result and the single precision result as the "computed" result, then the *number of* (*extra) accurate digits* can be defined as the overall average of $-\log_{10}\left|\frac{result_{true} - result_{computed}}{result_{true}}\right|$ Heyman and Reeves (1989) and Heyman and O'Leary (1995), computed this statistic for a set of test problems when computing the mean first passage time matrix. Heyman and Reeves (1989) considered four different procedures (state-reduction, Gaussian elimination, and two closed form matrix solutions) while in Heyman and O'Leary (1995) an *UL* factorization with normalisation related to a state reduction procedure was used. However in both of these papers their results were displayed in figures and no actual numerical results were tabulated.

As done for the stationary distributions, we also compute the single precision version of $M$ for each algorithm and comparing directly with the rounded double precision version and then calculating the average number of accurate decimal places obtained by the single precision version.

**7.3 The Group inverse**

A single measure for the accuracy of the group inverse is more problematic as there are three conditions that the group inverse must satisfy and a direct computation of the group inverse using matrix inversions is prone to multiple errors. Further we do not have any exact results published in the literature for group inverses of $I - P$.

The three conditions for parameterization of a generalized matrix inverse have been considered and three different statistics have been introduced to measure the accuracy of the computations. Let $A_G = I - (I - P)A^{\#}$, $B_G = I - A^{\#}(I - P)$. As in (2.6) and (2.8), let $\boldsymbol{\alpha} = A_G\boldsymbol{e}$, $\boldsymbol{\beta}^T = \boldsymbol{\pi}^T B_G$ and $\gamma = \boldsymbol{\beta}^T A^{\#}\boldsymbol{\alpha} - 1$.

From the parametrization of the group inverse, (Theorem 2, conditions (2.9) and (2.14)), these parameters are $\boldsymbol{\alpha} = \boldsymbol{e}$, $\boldsymbol{\beta} = \boldsymbol{\pi}$ and $\gamma = -1$.

Thus we compute, based on the single and double precision results, the following statistics: *MINDELTA* $\boldsymbol{\alpha} = \min_{1\le i\le m}|\alpha_i - 1|$, *MAXDELTA* $\boldsymbol{\alpha} = \max_{1\le i\le m}|\alpha_i - 1|$, *RELDELTA* $\boldsymbol{\alpha} = \sum_{i=1}^{m}|\alpha_i - 1|$, *MINDELTA* $\boldsymbol{\beta} = \min_{1\le i\le m}|\beta_i - \pi_i|$, *MAXDELTA* $\boldsymbol{\beta} = \max_{1\le i\le m}|\beta_i - \pi_i|$, *RELDELTA* $\boldsymbol{\beta} = \sum_{i=1}^{m}|\beta_i - \pi_i|$ and *DELTA* $\gamma = |\boldsymbol{\beta}A^{\#}\boldsymbol{\alpha}|$. All of these statistics should be close to zero.

In the paper of Heyman and O'Leary (1995), a procedure to compare the accuracy of the group inverse, based on the similar procedure as used for finding the average number of (extra) accurate digits for computing the mean first passage times, was used. We use this to find the *average number of (extra) accurate digits* in the computation of the group inverse in our perturbation procedures, using the technique we used for the mean first passage time matrix.

Further, as was done for the previous two properties, we also compute the single precision version of $A^{\#}$ for each algorithm and then compare this result directly with the appropriately rounded double precision elements to calculate the *average number of accurate decimal places obtained by the single precision version*.

We now illustrate the calculations performed for the following example.

**7.4 Test Example**: The transition matrix $P$ below appears in p.199 of Kemeny and Snell (1960) and in Sheskin (1985). The (1,1) entry was changed (as used in Hunter (1991)) to ensure that the matrix $P$ is in fact a stochastic matrix.

$$P = \begin{bmatrix} 0.831 & 0.033 & 0.013 & 0.028 & 0.095 \\ 0.046 & 0.788 & 0.016 & 0.038 & 0.112 \\ 0.038 & 0.034 & 0.785 & 0.036 & 0.107 \\ 0.054 & 0.045 & 0.017 & 0.728 & 0.156 \\ 0.082 & 0.065 & 0.023 & 0.071 & 0.759 \end{bmatrix}.$$

### 7.4.1 MatLab errors

For our numerical computation we use MatLab software, in particular the 64-bit version R2015b on a MacBook Air computer. In interpreting the errors generated by Matlab care must be taken. (See Moler (1996) and Moler (2004).) Under double precision (single precision) the spacing between 1 and the next largest number is the 2.2204e-16 (1.1921e-07), the eps. This floating point accuracy also depends on the size of the number and for numbers smaller than 1 (typically say probabilities) the spacing will be smaller. For example, observe that $a = 0.1 + 0.2 - 0.3$ yields a MatLab value of 5.5511e-17. This error is due to the machine precision which, of course, we cannot eliminate.

### 7.4.2 Stationary distribution

No exact results for the stationary distribution of the M. C. with transition matrix *P* appears in the literature. The stationary probability vector $\boldsymbol{\pi}^T$ using the *GTH* algorithm is, under double precision to 15 decimal places is given as $\pi^T(GTHD)$ = (0.270457577293538, 0.184235456501417, 0.076135265451860, 0.147597142335324, 0.321574558417861).

For all algorithms, under double precision, the calculations leading to the corresponding stationary probability distributions yield the same results as computed by the *GTH* algorithm, when rounded to 14 decimal places. (Actually the average number of accurate digits for the stationary distribution for each algorithm ranges from 14.6 to 14.8, as displayed in Table 1.)

The single precision version the stationary probability distribution given by *GTH* is same as the double precision version, when rounded to 6 decimal places. Further, the single precision versions of the stationary distributions derived using each algorithm also give the double precision *GTH* version when rounded to 6 decimal places for *AL*1 and *AL*2 but 7 decimal paces for A*L3*, *AL4A*, *AL4B* and *AL4C*. The average number of decimal places ranges from 6.6 to 7.8 with *AL4B* and *AL4C* both attaining 7.8.

As can been seen from Table 1, all algorithms yield very small errors for *MINRE(D)*, *MAXRE(D)* and *RELE(D)*, with perhaps *AL*2 and *AL4C* slightly inferior.

The relevant accuracy statistics show that, to 15 decimal places, we cannot detect any significant differences between *MINE, MAXE* and *RELE* for the pairs *(GTHD, D)* although *AL*2 and, especially, *AL4C* are marginally inferior. Similarly, under single precision, for the pairs *(GTHD, S), AL*1 and *AL*2 are marginally inferior to the other algorithms for *MINE, MAXE* and *RELE*.

It is difficult to identify and recommend one specific algorithm that performs better than others, although note that *AL4A* and *AL4B* (which are minor variants with different initial conditions) each perform consistently well across all error categories with mainly smaller errors than the other algorithms.

TABLE 1: Errors for stationary distributions under single and double precision.

| | *GTH* | *AL1* | *AL2* | *AL3* | *AL4A* | *AL4B* | *AL4C* |
|---|---|---|---|---|---|---|---|
| *MINRE(S)* | 1.2943e-09 | 5.0298e-09 | 0 | 2.5147e-09 | 3.6394e-11 | 1.4715e-10 | 2.9227e-11 |
| *MAXRE(S)* | 1.5001e-08 | 4.0083e-08 | 1.4901e-08 | 1.6366e-08 | 1.8086e-08 | 1.8200e-08 | 1.8069e-08 |
| *RELE(S)* | 3.9539e-08 | 8.3314e-08 | 1.4901e-08 | 4.0149e-08 | 3.4222e-08 | 3.4867e-08 | 3.4206e-08 |
| *Av # d.p.'s for π(S)* | | 7.2 | 6.6 | 7.4 | 7.6 | 2.3188e-08 | 2.2792e-08 |
| *MINE(GTHS, S)* | | 5.1810e-09 | 1.3499e-08 | 8.5802e-09 | 6.8943e-09 | 6.9936e-09 | 6.8804e-09 |
| *MAXE(GTHS, S)* | | 8.4687e-08 | 7.4605e-08 | 2.7866e-08 | 2.2766e-08 | | |
| *RELE(GTHS, S)* | | 1.7997e-07 | 1.6876e-07 | 7.2892e-08 | 6.7103e-08 | 6.7751e-08 | 6.7120e-08 |
| *MINE(S, D)* | 2.4537e-09 | 7.6347e-09 | 2.8958e-09 | 1.4306e-09 | 2.8898e-10 | 4.0312e-10 | 2.7164e-10 |
| *MAXE(S, D)* | 2.3830e-08 | 8.8179e-08 | 7.1114e-08 | 2.4374e-08 | 1.9275e-08 | 1.9697e-08 | 1.9300e-08 |
| *RELE(S, D)* | 5.4643e-08 | 1.7636e-07 | 1.40576e-07 | 5.1610e-08 | 3.8667e-08 | 3.9511e-08 | 3.8624e-08 |
| *MINE(GTHD, D)* | | 0 | 0 | 0 | 0 | 4.1633e-17 | 0 |
| *MAXE(GTHD, D)* | | 2.2204e-16 | 4.9960e-16 | 1.1102e-16 | 1.6653e-16 | 5.5511e-17 | 4.9960e-16 |
| *RELE(GTHD, D)* | | 6.3838e-16 | 9.7145e-16 | 2.9143e-16 | 3.0531e-16 | 2.6368e-16 | 1.1935e-15 |
| *MINRE(D )* | 0 | 0 | 0 | 0 | 0 | 0 | 0 |
| *MAXRE(D)* | 5.5511e-17 | 5.5511e-17 | 1.1102e-16 | 5.5511e-17 | 5.5511e-17 | 0 | 1.1102e-16 |
| *RELE(D)* | 5.5511e-17 | 1.6653e-16 | 2.2204e-16 | 8.3267e-17 | 6.9389e-17 | 0 | 2.7756e-16 |
| *Av # d.p.'s for π(D)* | 15 | 14.8 | 14.8 | 14.8 | 14.6 | 14.8 | 14.8 |

### 7.4.3. Mean first passage time matrix *M*:

For each of the different algorithms the mean first passage times are given (to 12 decimal places) as follows:

$$M = \begin{bmatrix} 3.697437542727 & 22.374164571709 & 57.756742192108 & 23.278850538432 & 9.598732858601 \\ 17.032615490720 & 5.427836850679 & 56.864516889123 & 22.100075015307 & 8.844407674651 \\ 17.667201055109 & 22.106202543394 & 13.134517809389 & 22.292628444747 & 9.020416501550 \\ 16.341175493452 & 21.005100548563 & 56.552837505099 & 6.775198924435 & 7.609106618566 \\ 15.243523199997 & 20.060109096789 & 55.798746557709 & 20.158095744297 & 3.109698742711 \end{bmatrix}.$$

TABLE 2: Errors for mean first passage times under single and double precision.

| | *AL1* | *AL2* | *AL3* | *AL4A* | *AL4B* | *AL4C* |
|---|---|---|---|---|---|---|
| *MINRESM(S)* | 0 | 0 | 0 | 0 | 0 | 0 |
| *MAXRESM(S)* | 1.7285e-06 | 3.8147e-06 | 3.7104e-06 | 4.4256e-06 | 2.7269e-06 | 4.3064e-06 |
| *RESM(S)* | 1.6227e-05 | 1.9968e-05 | 1.9461e-05 | 2.1696e-05 | 1.8865e-05 | 1.9073e-05 |
| *Accurate d.p.'s for M(S)* | 5.04 | 5.00 | 5.20 | 5.56 | 4.88 | 5.28 |
| *MINEM(S, D)* | 4.7146e-08 | 2.4245e-08 | 4.2147e-08 | 3.04986e-09 | 1.1877e-08 | 8.0028e-09 |
| *MAXEM(S, D)* | 1.3477e-05 | 1.8601e-05 | 3.7823e-06 | 4.5595e-06 | 8.0312e-06 | 4.4580e-06 |
| *REM(S, D)* | 8.0903e-05 | 1.1128e-04 | 2.5273e-05 | 2.3256e-05 | 4.2810e-05 | 2.4030e-05 |
| *MINRESM(D)* | 0 | 0 | 2.7756e-17 | 5.551e-17 | 0 | 1.1102e-16 |
| *MAXRESM(D)* | 1.3378e-14 | 7.2164e-15 | 7.8826e-15 | 6.9389e-15 | 6.7168e-15 | 8.0214e-15 |
| *RESM(D)* | 6.5808e-14 | 4.9655e-14 | 5.5261e-14 | 3.6221e-14 | 3.8386e-14 | 4.8128e-14 |
| *Accurate Digits* | 7.0475 | 7.1011 | 7.5035 | 7.8174 | 7.5678 | 7.6773 |

Calculations in single precision in general do poorly and we typically achieve much less accuracy than the expected six decimal places. (see Table 2.) *AL*1 gives the most accurate results with *AL*4*A* the worst.

By comparing the single precision results for each of the algorithms against rounding the terms for *M* above we achieve an average of 5.56 decimal places for *AL4A* with smaller averages for the other algorithms (with *AL4*B being the worst at an average of 4.88 decimal places).

In terms of double precision, all *AL* achieve relatively small *MINRESM(D)* errors. *AL*1 has the largest *MAXRESM(D)* and *RESM(D)* errors. Even though *AL4B* gives the smallest *MAXRESM(D)*, *AL4A* gives the smallest *MINEM(S, D), REM(S, D), RESM(D)* errors as well as the largest number of accurate digits.

Overall *AL4A*, except for the single precision results, appears to perform consistently well and is recommended for calculating the mean first passage times.

**7.4.4. Group inverse $A^{\#}$:**
The group inverse for $I - P$ is given (to 13 decimal places) is as follows: $A^{\#} =$

$$\begin{bmatrix} 3.1905741863522 & -0.9375239582265 & -0.4087732024356 & -0.6983862380226 & -1.1458907876676 \\ -1.4160257342402 & 3.1845904654802 & -0.3408433921500 & -0.5244023393545 & -0.9033189997355 \\ -1.5876542085704 & -0.8881558516147 & 3.9885516959952 & -0.5528226752867 & -0.9599189605234 \\ -1.2290205477352 & -0.6852938229425 & -0.3171135995114 & 2.7375055783011 & -0.5060776081121 \\ -0.9321521677369 & -0.5111928914350 & -0.2597006850570 & -0.2377717484790 & 1.9408174927079 \end{bmatrix}.$$

All the algorithms give the above expression to the requisite decimal places when the computations are carried out in double precision.

When the terms of $A^{\#}$ above are rounded to 5 decimal places we obtain the group inverse expression for all the algorithms when the calculations are carried out under single precision (with one entry for *AL*1 only at 4 decimal places). However, for every algorithm there are some entries that can be expressed to seven decimal places (and even one to 8 d.ps for *AL4B*). In Table 3 we give figures for the average number of accurate decimal places with *AL4B* having the highest average of 6.28, under single precision.

Table 3 gives the accuracy of the $\alpha, \beta$ and $\gamma$ parameters for the parametric form for the group generalized inverse $A^{\#}$ derived using each algorithm.

TABLE 3: Errors for the group inverse under single and double precision.

| | *AL*1 | *AL2* | *AL3* | *AL4A* | *AL4B* | *AL4C* |
|---|---|---|---|---|---|---|
| *MINDELTA* $\alpha(S)$ | 0 | 0 | 0 | 0 | 0 | 0 |
| *MAXDELTA* $\alpha(S)$ | 1.1921e-07 | 1.1921e-07 | 0 | 1.1921e-07 | 1.1921e-07 | 0 |
| *RELDELTA* $\alpha(S)$ | 1.7881e-07 | 1.1921e-07 | 0 | 1.7881e-07 | 1.7881e-07 | 0 |
| *MINDELTA* $\beta(S)$ | 8.0763e-09 | 0 | 1.4652e-09 | 2.7727e-09 | 1.4715e-10 | 2.7424e-09 |
| *MAXDELTA* $\beta(S)$ | 1.2473e-08 | 0 | 1.2865e-08 | 7.7656e-09 | 2.1615e-08 | 7.7911e-09 |
| *RELDELTA* $\beta(S)$ | 5.1601e-08 | 0 | 2.5248e-08 | 2.6771e-08 | 3.7371e-08 | 2.6756e-08 |
| *DELTA* $\gamma(S)$ | 2.2352e-08 | 1.4901e-08 | 8.9407e-08 | 6.7055e-08 | 3.7253e-08 | 8.9407e-08 |
| *Av # accurate d.p.'s* $A^{\#}(S)$ | 5.84 | 5.88 | 6.08 | 6.16 | 6.28 | 6.16 |
| *MINDELTA* $\alpha(D)$ | 0 | 0 | 1.1102e-16 | 0 | 0 | 0 |
| *MAXDELTA* $\alpha(D)$ | 2.2204e-16 | 2.2204e-16 | 6.6613e-16 | 2.2204e-16 | 2.2204e-16 | 4.4409e-16 |
| *RELDELTA* $\alpha(D)$ | 3.3307e-16 | 3.3307e-16 | 1.4433e-15 | 5.5511e-16 | 7.7716e-16 | 8.8818e-16 |
| *MINDELTA* $\beta(D)$ | 0 | 0 | 0 | 0 | 0 | 0 |
| *MAXDELTA* $\beta(D)$ | 5.5511e-17 | 2.7756e-17 | 5.5511e-17 | 5.5511e-17 | 5.5511e-17 | 5.5511e-17 |
| *RELDELTA* $\beta(D)$ | 9.7145e-17 | 2.7756e-17 | 5.5511e-17 | 8.3267e-17 | 6.9389e-17 | 8.3267e-17 |
| *DELTA* $\gamma(D)$ | 3.7470e-16 | 4.1633e-17 | 3.4694e-16 | 9.7145e-17 | 2.2204e-16 | 4.1633e-17 |
| *Accurate digits* | 6.7495 | 6.8348 | 7.0582 | 7.1085 | 7.2211 | 7.1458 |

Under single precision, *AL*3 and *AL4C* both give *DELTA* $\alpha$ explicity 0. *AL*2 gives *DELTA* $\beta$ explicitly 0 with *AL*1 the least accurate. *AL2* gives the smallest error for $\gamma$.

Under double precision, all *AL*, except for *AL3* have zero *MINDELTA* $\alpha$ errors. *AL*1 and *AL*2 give the smallest overall *RELDELTA* $\alpha$ error. All *AL* have zero MIN*DELTA* $\beta$ errors with *AL*2 having the smallest *MAXDELTA* $\beta$ and *AL*1 the largest. For *DELTA* $\gamma$, *AL*2 and *AL4C* have the smallest error with *AL*1 the largest. Note however that for all algorithms under double precision we get accuracy to at least 15 d.p.'s for all the parameters $\alpha$, $\beta$ and $\gamma$.

While it is difficult to give a universal recommendation based on the above observations, overall, *AL*2 appears to give consistently accurate results for all the parameters although *AL4B* gives the most accurate number of digits, under both single and double precision.

An interesting observation is that no one particular algorithm for computing the key Markov chain properties has emerged to dominate the accuracy of all the different procedures.

In a sequel paper, when we compare not only the perturbation procedures but alternative computational techniques for the key properties of irreducible M. C.'s, we may be able to gain a better impression as to whether perturbation procedures may in fact prove to be suitable alternatives.

**Acknowlegements**
This research was initiated following discussions with Professor Stephen Kirkland when the author visited him at the National University of Ireland, Maynooth. Steve suggested the procedure described by Algorithm 2 that arose from ideas in the book he wrote with the late Miki Neumann, (Kirkland, S. J., Neumann, M. (2013)). The author wishes to express his thanks for his contribution, friendship and the hospitality extended during his visit in May 2012.
The author would like to express his thanks to Ms Diane Park who coded the algorithms and implemented some of the MatLab calculations leading to the results given in this paper. This was part of a Summer Research project that lead on to her completing a BSc(Hons) degree at Auckland University of Technology with First Class Honours in Mathematics in 2014.
The author also wishes to express his appreciation to the referee. In particular, the referee pointed out a gap in the proof of Theorem 10 that had gone undetected in an earlier paper by the author and provided a range of very useful comments and additional references.

**References**

Ben-Ari, I., Neumann M. (2012). Probabilistic approach to Perron root, the group inverse, and applications, *Linear Multilinear A*, 60, 39-63.

Berman, A., Plemmons R.J. (1979). *Nonnegative Matrices in the Mathematical Sciences*. Academic Press, New York. Reprinted (1994), SIAM, Philadelphia.

Feller, W. (1950), *An Introduction to Probability Theory and its Applications, Volume* 1, Wiley, New York, 2nd Edition (1957), 3rd Edition (1967).

Grassman, W. K., Taksar, M. I., Heyman, D. P. (1985). Regenerative analysis and steady state distributions for Markov chains. *Oper Res*, **33**, 1107-1116.

Funderlic, R., Plemmons, R. (1986). Updating LU factorizations for computing stationary distributions. *SIAM J Alg Disc Meth*, **7**, 30-42.

Heyman, D. P., O'Leary, D. P. (1995). What is Fundamental for Markov Chains: First Passage Times, Fundamental Matrices, and Group Generalized Inverses. *Computations with Markov chains: Proceedings of the 2nd International Workshop on the Numerical Solution of Markov Chains*. Stewart, W. J. (Ed.). pp151-161 Kluwer Academic Publishers, Dordrecht.

Harville, D. A. (1997). *Matrix Algebra from a Statistician's Perspective*. Springer-Verlag, New York.

Heyman D. .P., Reeves A. (1989). Numerical solutions of Linear Equations Arising in Markov chain models. *ORSA J Comput*, 1, 52-60.

Hunter, J. J. (1969). On the moments of Markov renewal processes. *Adv Appl Probab*, **1**, 188-210.

Hunter, J. J. (1982). Generalized inverses and their application to applied probability problems. *Linear Algebra Appl.* **45**, 157-198.

Hunter, J. J. (1986). Stationary distributions of perturbed Markov chains. *Linear Algebra Appl.* **82**, 201-214.

Hunter, J. J. (1988). Characterisations of generalized inverses associated with Markovian kernels. *Linear Algebra Appl.* **102**, 121-142.

Hunter, J. J. (1990). Parametric forms for generalized inverses of Markovian kernels and their applications. *Linear Algebra Appl.* **127**, 71-84.

Hunter, J. J. (1991). The computation of stationary distributions of Markov chains through perturbations. *J Appl Math Stoch Anal*, **4**, 29-46.

Hunter, J. J. (1992). Stationary distributions and mean first passage times in Markov chains using generalized inverses. *Asia Pac J Oper Res*, **9**, 145-153.

Hunter, J. J. (2006). Mixing times with applications to perturbed Markov chains. *Linear Algebra Appl*, **417**, 108-123.

Hunter, J. J. (2007). Simple procedure for finding mean first passage times in Markov chains. *Asia Pac J Oper Res*, **24**, 813-829.

Hunter, J. J. (2014). Generalized inverses of Markovian kernels in terms of properties of the Markov chain. *Linear Algebra Appl,* **447**, 38-55.

Hunter, J. J. (2016). Accurate calculations of stationary distributions and mean first passage times in Markov renewal processes and Markov chains. *Special Matrices*, **4**, 151-175.

Kemeny, J. G., Snell, J. L. (1960*). Finite Markov Chains.* Van Nostrand, New York.

Kohlas, J. (1986). Numerical computation of mean first passage times and absorption probabilities in Markov and semi-Markov models. *Zeit Fur Oper Res*, **30**, 197-207.

Kirkland, S. J., Neumann, M. (2013). *Group Inverses of M-matrices and their applications*. CRC Press. Taylor & Francis Group, Boca Raton, Florida, U.S.A.

Langville, A., Meyer, C. D. Jr. (2006). Updating Markov chains with an eye on Google's PageRank. *SIAM J Matrix Anal Applicat*. **27**, 968-987.

Meyer, C. D. Jr. (1975). The role of the group generalized inverse in the theory of finite Markov chains. *SIAM Rev,* **17**, 443-464.

Meyer, C. D. Jr. (1980). The condition of a finite Markov chain and perturbation bounds for the limiting probabilities. *SIAM J Alg Disc Meth*, **1**, 273-283.

Meyer, C. D. Jr., Shoaf, J. (1980). Updating finite Markov chains by using group matrix inversion. *J Stat Comput Sim*, **11**, 163-181.

Moler, C. (1996). Floating Points, *MATLAB News and Notes, Fall*, 1996. A PDF version is available at *http://www.mathworks.com/company/newsletters/news_notes/pdf/Fall96Cleve.pdf*

Moler, C. (2004). *Numerical Computing with MATLAB*, S.I.A.M. A PDF version is available on the MathWorks website at http://www.mathworks.com/moler/.

Paige, C.C., Styan, G.P. H., Wachter, P. G., (1975). Computation of the stationary distribution of a Markov chain. *J Stat Comput Sim*, **4**, 173-186.

Seneta, E. (1991). Sensitivity analysis, ergodicity coefficients, and rank-one updates for finite Markov chains. In: Stewart, W. J. (Ed) *Numerical Solution of Markov Chains, Probability: Pure and Applied,* Marcel-Dekker, New York, 121-129.

Sherman, J., Morrison, W. J. (1949). Adjustment of an Inverse Matrix Corresponding to Changes in the Elements of a Given Column or a Given Row of the Original Matrix (abstract). *Ann Math Stat,* **20**: 621.

Sheskin, T.J. (1985). A Markov chain partitioning algorithm for computing steady state probabilities. *Oper Res,* **33**, 228-235.